\documentclass[12pt]{article} % use larger type; default would be 10pt

\usepackage[utf8]{inputenc} % set input encoding (not needed with XeLaTeX)

\usepackage{geometry} % to change the page dimensions
\usepackage{graphicx} % support the \includegraphics command and options

\usepackage{booktabs} % for much better looking tables
\usepackage{array} % for better arrays (eg matrices) in maths
\usepackage{paralist} % very flexible & customisable lists (eg. enumerate/itemize, etc.)
\usepackage{verbatim} % adds environment for commenting out blocks of text & for better verbatim
\usepackage{subfig} % make it possible to include more than one captioned figure/table in a single float
\usepackage{fancyhdr} % This should be set AFTER setting up the page geometry
\usepackage{sectsty}
\allsectionsfont{\sffamily\mdseries\upshape} % (See the fntguide.pdf for font help)
\usepackage[nottoc,notlof,notlot]{tocbibind} % Put the bibliography in the ToC
\usepackage[titles,subfigure]{tocloft} % Alter the style of the Table of Contents

\usepackage{url} % For putting url-s inside the bititems
\usepackage{amsthm, amssymb,amsmath} % I am installing this package in order to be able to write theorems with proofs
\title{Hamiltonian Subgraphs of Order Seven in $srg(n,k,1,2)$}
\author{Reimbay Reimbayev}
\date{} % Activate to display a given date or no date (if empty),
\begin{document}
\maketitle

\begin{abstract}
Strongly regular graphs are highly symmetrical and can be described fully with just a few parameters, yet the existence of many of them is still under the question. In this paper, we continue the study of the famuly of strongly regular graphs with parameters $\lambda =1$ and $\mu =2$ and establish all of their possible Hamiltonian subgraphs of order seven. By doing so we establish the lower and upper bounds for number of 7-gons, or 7-cycles, in such graphs.
\end{abstract}

%\section{Introduction}
\bigskip
We will skip a long and lengthy introduction into the problem \cite{Conway} and motivation of this work as this is merely a continuation of the previous work that we have posted on ArXiv. But instead, we will briefly summarize what we have done so far. In our previous papers, we have extensively studied the structure of the strongly regular graphs with parameters $\lambda=1$ and $\mu=2$ \cite{ReiLowerBound, ReiAllSix}. In particular, we have looked at the possible subgraphs of order six and their numbers in such graphs, should they exist. Here, we continue our previous work and move on subgraphs of order seven but only the ones that contain Hamiltonian cycles, or in short, the Hamiltonian subgraphs of order seven. By doing so, we will also establish bounds for heptagons, $C_7$, in addition to other polygons (cycles) of lower order that we have found previously.

\begin{figure}
	\includegraphics[width=1.0\textwidth]{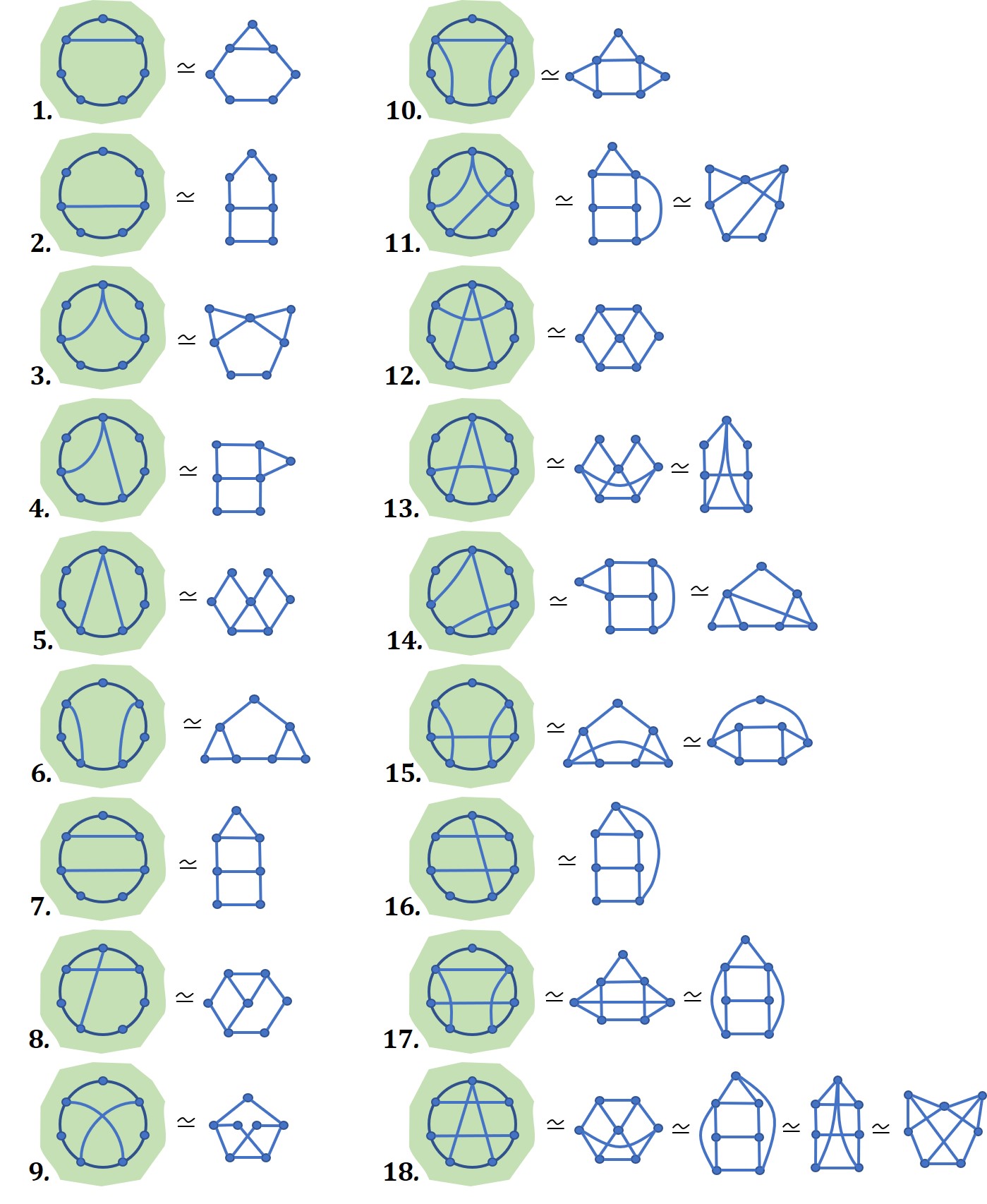}
		\centering
		\caption{All possible Hamiltonian subgraphs, except one for a 7-cycle, in $srg(n,k,1,2)$.}
		\label{mainFigure}
\end{figure}

%\section{Main Results}
\bigskip
Strongly regular graphs with parameters $\lambda=1$ and $\mu=2$, for simplicity henceforth call it a graph $G$, can have 19 possible Hamiltonian subgraphs of order seven (Figure \ref{mainFigure}). The cycle of order seven is not depicted in the figure. Obviously, for the graphs of smaller order (smaller $n$) not all the subgraphs would exist. The subgraphs are obtained by straight-forward exhaustive search of all such graphs that do not brake strong-regularity conditions for $G$.

Let us denote $H_i$ the graphs of type $i$ from the figure \ref{mainFigure} and by $h_i$ the number of such graphs in $G$. Separately, we denote $H_0$ the cycles of order seven, $C_7$, and $h_0$ -  the number of such cycles in $G$. We will also extensively use the notation from the previous study regarding six-vertex subgraphs $N_i$ and their number $n_i$ given in the paper \cite{ReiAllSix}. We assume they are already known. Here we derive only the values for $h_i$.

Some of the values for $h_i$-s can be found directly. For example, $h_{10}$. For that consider $N_3$. On one of the two sides of its quadrilateral with no triangles attached , we can recover a triangle in a unique way. Thus,
\[ h_{10}=2n_3.\]

Similarly, $H_{17}$ can be obtained from $N_1$ in a unique way, by recovering a triangle on one of its three (free of triangles) sides. Thus,
\[h_{17}=3n_1=\frac{1}{4}nk(k-2)-n_3.\]

Now $H_{14}$ can be obtained from $N_4$ uniquely in two ways:
\[h_{14}=2n_4=4n_3.\]

Next, consider $N_5$. The subgraph has two triangles. Each vertex of degree two on one triangle has exactly two common neighbors with any two vertices of degree two of the other triangle. Thus,
\[h_6=8n_5=nk(k-2)(k-4)-8n_3.\]

From $N_9$, we can recover $H_4$. Thus,
\[h_4=nk(k-2)(k-4)-4n_3.\]

Consider $N_8$. Recover a triangle at one of its edges that are not incident with its existing triangle (or belong to it). As there are only two such edges, there are two ways of doing that. We obtain:
\[2n_8=2h_6+4h_{15}+h_{14}.\]
Thus, \[4h_{15}=2n_8-2h_6-h_{14}=8n_3.\]
Or, \[h_{15}=2n_3.\]

To move further, we have to assign one of the $h_i$-s as a free variable. The choice is fairly arbitrary. Here we have chosen $h_{11}$.

Consider $N_4$. Recover a triangle on its edge that is incident to the existing triangle of $N_4$, but only on the ones that are on its sides (leaving the one in the middle). Notice, that the middle edge is topologically different from the two on the sides, as it is incident to vertices of degrees both equal three. Thus,
\[2n_4=h_{11}+4h_{18}.\]
Or,
\[h_{18}=n_3-\frac{h_{11}}{4}.\]

In the same way we recover triangles on $N_8$ and obtain:
\[2n_8=2h_3+h_{11}.\]
Or, \[h_3=nk(k-2)(k-4)-2n_3-\frac{h_{11}}{2}.\]

Consider again $N_4$. It has a pair of vertices of degree two which must have exactly one more neighbor out of $N_4$. Recover it. We have,
\[n_4=h_{16}+4h_{18}.\]
Or, \[h_{16}=-2n_3+h_{11}.\]

As a side note, from the values of $h_{16}$ and $h_{18}$, both nonnegative integers, we can give an estimates for $h_{11}$:
\[4n_3 \geq h_{11} \geq 2n_3.\]

Once more, consider $N_4$. This time we are looking for the second common neighbor between a vertex of degree two (there are two of them but we take one) and the vertex of the triangle that doesn't belong to the same equilateral as the first vertex. We have,
\[2n_4=2h_{13}+4h_{18}.\]
Or,\[h_{13}=\frac{h_{11}}{2}.\]

Consider, $M_{19}$, or basically two triangles of $G$ that share a common vertex. The value of $m_{19}$ is also taken from \cite{ReiAllSix}. We can build it up to $h_{12}$ and by doing that we can also obtain $H_{18}$. Thus,
\[2m_{19}=h_{12}+h_{18}.\]
Or, \[h_{12}=\frac{1}{4}nk(k-2)-n_3+\frac{h_{11}}{4}.\]

Consider $N_9$. The subgraph has exactly two pairs of vertices that are neither adjacent nor have a common neighbor in $N_9$. To a chosen pair we add one of their common neighbor from $G$. We obtain:
\[4n_9=h_9+2h_{16}+2h_{18}.\]
Or, \[h_9=nk(k-2)(k-4)-2n_3-\frac{3}{2}h_{11}.\]

Consider $N_8$. For the vertex of degree two of the triangle, there are exactly two vertices of $N_8$ that are not adjacent to it but have a common neighbor (one of the vertices of the triangle). Recover the second common neighbor from $G$. We get:
\[2n_8=h_8+2h_{16}+4h_{15}.\]
Or,\[h_8=2nk(k-2)(k-4)-8n_3-2h_{11}.\]

Consider $N_9$. Recover a triangle on one of the two edges that are incident to vertices of degree of two from both sides. We have:
\[2n_9=h_7+2h_{16}+2h_{18}.\]
Or, \[h_7=\frac{1}{2}nk(k-2)(k-4)-\frac{3}{2}h_{11}.\]

Consider $N_{17}$. The vertex of degree two of the triangle of $N_{17}$ (the `rooftop') and its leaf has exactly one common neighbor outside of $N_{17}$. Recover it. We have:
\[n_{17}=2h_5+2h_{13}.\]
Or, \[h_5=\frac{1}{2}nk(k-2)(k-4)-\frac{h_{11}}{2}.\]

Consider $N_{11}$. From the set of vertices $G \backslash N_{11}$ recover a common neighbor for the leaf and one of the vertices of degree two of the quadrilateral that is adjacent to the vertex of degree three in $N_{11}$. We have,
\[4n_{11}=2h_2+2h_7+2h_5+h_8+h_{11}.\]
Or, \[h_2=nk(k-2)(k-4)(k-8)+12n_3+\frac{5}{2}h_{11}.\]

Consider a hexagon, $N_{12}$. Recover a triangle on one of its sides. We have:
\[6n_{12}=h_1+h_8+2h_{12}.\]
Or, \[h_1=\frac{1}{2}nk(k-2)(2k^2-25k+68)+16n_3+\frac{3}{2}h_{11}.\]

Finally, for $h_0$, the number of seven-cycles in $G$, or heptagons, which are not depicted in the figure, we have a relation:
\[ 2n_{35}=7h_0+2h_1+2h_2+h_3+h_4+h_5.\]
Thus, \[h_0=\frac{1}{14}nk(k-2)(k-4)(2k^2-30k+133)-10n_3-h_{11}.\]

Notice, that given the values for $n_3$ and $h_{11}$ cannot be negative, we can see what should be the upper bound for number of heptagons in $G$. It is not hard to derive the lower bound as well but we would rather leave it as it is as we strongly believe that the upper bound is indeed its exact value. To prove that still requires some meticulous work so for now we leave it as a conjecture. 

\bigskip
Let us summarize what we have so far. Below are all the values for number of Hamiltonian subgraphs in $srg(n,k,1,2)$. They are given in terms of `free variables' $n_3$ and $h_{11}$. Of course, it is possible using the relation between $n$ and $k$ to get rid of $n$ altogether but then the formulas become very cumbersome.

\begingroup
\allowdisplaybreaks
 \begin{flalign*}
 h_0=&\frac{1}{14}nk(k-2)(k-4)(2k^2-30k+133)-10n_3-h_{11};\\
 h_1=&\frac{1}{2}nk(k-2)(2k^2-25k+68)+16n_3+\frac{3}{2}h_{11};\\
 h_2=&nk(k-2)(k-4)(k-8)+12n_3+\frac{5}{2}h_{11};\\
 h_3=&nk(k-2)(k-4)-2n_3-\frac{h_{11}}{2};\\
 h_4=&nk(k-2)(k-4)-4n_3;\\
 h_5=&\frac{1}{2}nk(k-2)(k-4)-\frac{h_{11}}{2};\\
 h_6=&nk(k-2)(k-4)-8n_3;\\
 h_7=&\frac{1}{2}nk(k-2)(k-4)-\frac{3}{2}h_{11};\\
 h_8=&2nk(k-2)(k-4)-8n_3-2h_{11};\\
 h_9=&nk(k-2)(k-4)-2n_3-\frac{3}{2}h_{11};\\
 h_{10}=&2n_3;\\
 h_{11}=&h_{11};\\
 h_{12}=&\frac{1}{4}nk(k-2)-n_3+\frac{h_{11}}{4};\\
 h_{13}=&\frac{h_{11}}{2};\\
 h_{14}=&4n_3;\\
 h_{15}=&2n_3;\\
 h_{16}=&h_{11}-2n_3;\\
 h_{17}=&\frac{1}{4}nk(k-2)-n_3;\\
 h_{18}=&n_3-\frac{h_{11}}{4}.\\
\end{flalign*}
\endgroup

Finally, denote $p_i$ the number of i-gons in $G$. By i-gon here we mean cycle $C_i$ in $G$. Then, using the previously obtained results \cite{ReiLowerBound, ReiAllSix} for $p_i$-s up to six and adding one more for $p_7$, we have the following relations:
\begin{align*}
p_3 =&\frac{1}{6}nk;\\
p_4 =& \frac{1}{8}nk(k-2);\\
p_5=&\frac{1}{5}nk(k-2)(k-4);\\
p_6\geq&\frac{1}{12}nk(k-2)(2k^2-21k+53);\\
p_7\leq&\frac{1}{14}nk(k-2)(k-4)(2k^2-30k+133).
\end{align*}
 
We conjecture that the bounds for $p_6$ and $p_7$ are in fact their exact values.

%\subsection{A subsection}

\end{document}